\documentclass[reqno]{amsart}
\usepackage{hyperref,amsmath,amsthm,amssymb}
\usepackage[dvips]{graphicx}

\theoremstyle{plain}
\newtheorem{thm}{Theorem}
\newtheorem{lem}{Lemma}
\newtheorem{cor}{Corollary}
\newtheorem{ex}{Example}

\renewcommand{\Re}{\mathrm{Re}}
\renewcommand{\Im}{\mathrm{Im}}

\title{An extension of Nunokawa lemma and its example}

\author{Hitoshi Shiraishi}
\address{Hitoshi Shiraishi \newline
Department of Mathematics \newline
Kinki University \newline
Higashi-Osaka, Osaka 577-8502, Japan}
\email{step\_625@hotmail.com}

\author{Mamoru Nunokawa}
\address{Mamoru Nunokawa \newline
Emeritus Professor \newline
University of Gunma \newline
798-8, Hoshikuki, Chuou, Chiba 260-0808, Japan.}
\email{mamoru\_nuno@doctor.nifty.jp}

\date{}

\begin{document}
\maketitle

\begin{abstract}
For analytic functions $p(z)$ in the open unit disk $\mathbb{U}$ with $p(0)=1$,
Nunokwa has given a result which called Nunokawa lemma (Proc. Japan Acad., Ser. A {\bf 68} (1992)).
By studying Nunokawa lemma,
we obtain this expansion.
In this paper, we introduce this result and its example.
\end{abstract}

\

\section{Introduction}

\

Let $\mathbb{U}$ be defined by the open unit disk
$$
\mathbb{U}
= \{z \in \mathbb{C}:|z|<1\}.
$$

\

The basic tool in proving our results is the following lemma due to Miller and Mocanu \cite{m1ref2} (also \cite{d7ref2}).

\

\begin{lem} \label{jack} \quad
Let the function $w(z)$ be analytic in $\mathbb{U}$ with $w(0)=0$.
If $\left|w(z)\right|$ attains its maximum value on the circle $|z|=r$ at a point $z_{0}\in\mathbb{U}$,
then there exists a real number $m \geqq 1$ such that
$$
\frac{z_{0}w'(z_{0})}{w(z_{0})}
= m.
$$
\end{lem}

\

\section{Main result}

\

Applying Lemma \ref{jack},
we derive the following result.

\

\begin{thm} \label{p01thm1} \quad
Let $p(z)$ be analytic in $\mathbb{U}$ with $p(0)=1$ and suppose that there exists a point $z_0\in\mathbb{U}$ such that
$$
\Re(p(z))
> \alpha
\quad for \quad |z|<|z_0|
$$
and
$$
p(z_0)
= \alpha + \beta i
$$
for some real $\alpha$ and $\beta$, $0\leqq\alpha<1$ and $\beta\neq0$.

Then we have
$$
\Re\left( \frac{z_0 p'(z_0)}{p(z_0)} \right)
= -\frac{\alpha \beta k}{\alpha^2+\beta^2}
\leqq 0
$$
and
$$
\Im\left( \frac{z_0 p'(z_0)}{p(z_0)} \right)
= \frac{\beta^2 k}{\alpha^2+\beta^2}
$$
where
$$
k
\geqq \frac{1}{2}\left( \frac{\beta}{1-\alpha}+\frac{1-\alpha}{\beta} \right)
\geqq 1
\qquad(\beta>0)
$$
and
$$
k
\leqq \frac{1}{2}\left( \frac{\beta}{1-\alpha}+\frac{1-\alpha}{\beta} \right)
\leqq -1
\qquad(\beta<0).
$$
\end{thm}

\

\begin{proof} \quad
Let us define the function $q(z)$ by
$$
q(z)
= \frac{p(z)-\alpha}{1-\alpha}
\qquad(z\in\mathbb{U}).
$$

Clearly, $q(z)$ is analytic in $\mathbb{U}$ with $q(0)=1$ and $q(z_0)=\dfrac{\beta}{1-\alpha}i$ for a point $z_0$ such that $p(z_0)=\alpha+\beta i$.

Also, let us put
$$
w(z)
= \frac{1-q(z)}{1+q(z)}
\qquad(z\in\mathbb{U}).
$$

Then, we have that $w(z)$ is analytic in $|z|<|z_0|$, $w(0)=0$, $|w(z)|<1$ for $|z|<|z_0|$ and
$$
|w(z_0)|
= \left| \frac{(1-\alpha)^2-\beta^2-2(1-\alpha)\beta i}{(1-\alpha)^2+\beta^2} \right|
= 1.
$$

From Lemma \ref{jack},
we obtain
$$
\frac{z_0 w'(z_0)}{w(z_0)}
= \frac{-2 z_0 q'(z_0)}{1- \{ q(z_0) \}^2}
= \dfrac{-2 z_0 q'(z_0)}{1+ \left( \dfrac{\beta}{1-\alpha} \right)^2}
= m \geqq 1.
$$

This shows that
$$
-z_0 q'(z_0)
\geqq \frac{1}{2}\left( 1+ \left( \dfrac{\beta}{1-\alpha} \right)^2 \right)
$$
and $z_0 q'(z_0)$ is a negative real number.

From the fact that $z_0 q'(z_0)$ is a real number and $q(z_0)$ is a pure imaginary number,
we can put
$$
\frac{z_0q'(z_0)}{q(z_0)}
= ik
$$
where $k$ is a real number.

For the case $\beta>0$,
we have
\begin{align*}
k
&= \Im\left( \frac{z_0q'(z_0)}{q(z_0)} \right) \\
&= \Im\left( -z_0q'(z_0) \frac{1-\alpha}{\beta}i \right) \\
&\geqq \frac{1}{2}\left( 1+ \left( \dfrac{\beta}{1-\alpha} \right)^2 \right)\frac{1-\alpha}{\beta} \\
&= \frac{1}{2}\left( \frac{\beta}{1-\alpha}+\frac{1-\alpha}{\beta} \right)
\geqq 1
\end{align*}
and for the case $\beta<0$,
we get
\begin{align*}
k
&= \Im\left( \frac{z_0q'(z_0)}{q(z_0)} \right) \\
&= \Im\left( -z_0q'(z_0) \frac{1-\alpha}{\beta}i \right) \\
&\leqq \frac{1}{2}\left( 1+ \left( \dfrac{\beta}{1-\alpha} \right)^2 \right)\frac{1-\alpha}{\beta} \\
&= \frac{1}{2}\left( \frac{\beta}{1-\alpha}+\frac{1-\alpha}{\beta} \right)
\leqq -1.
\end{align*}

On the other hand,
let us consider
$$
\frac{z_0 q'(z_0)}{q(z_0)}
= \frac{z_0 p'(z_0)}{p(z_0)-\alpha}
= ik,
$$
then we have
$$
\frac{z_0p'(z_0)}{p(z_0)}
= \frac{p(z_0)-\alpha}{p(z_0)}ik
= -\frac{\alpha \beta k}{\alpha^2+\beta^2}+\frac{\beta^2 k}{\alpha^2+\beta^2}i.
$$

This completes our proof.
\end{proof}

\

Putting $\alpha=0$ in Theorem \ref{p01thm1},
we have Corollary \ref{p01cor1} \cite{ds2ref2}.

\

\begin{cor} \label{p01cor1} \quad
Let $p(z)$ be analytic in $\mathbb{U}$ with $p(0)=1$ and suppose that there exists a point $z_0\in\mathbb{U}$ such that
$$
\Re(p(z))
> 0
\quad for \quad |z|<|z_0|,
$$
$\Re(p(z_0))=0$ and $p(z_0)\neq0$.

Then we have
$$
\frac{z_0 p'(z_0)}{p(z_0)}
= ik
$$
where $k$ is a real and $k\geqq1$ for $\Im(p(z_0))>0$ and $k\leqq-1$ for $\Im(p(z_0))<0$.
\end{cor}

\

\section{Example of the theorem}

\

\begin{ex} \upshape \label{p01ex1} \quad
We consider the function $p(z)$ given by
$$
p(z)
= 1+(1-\alpha)(2z+z^2)
\qquad(z\in\mathbb{U})
$$
for some real $0\leqq\alpha<1$.

Then, $p(z)$ is analytic in $\mathbb{U}$ with $p(0)=1$.

Putting $z_0=-\dfrac{1}{2}\pm\dfrac{1}{2}i$,
it follows that
$$
\Re(p(z))
> \alpha
\quad for \quad |z|<|z_0|=\frac{1}{\sqrt[]{2}}
$$
and $\Re(p(z_0))=\alpha$.
For the case $z_0=-\dfrac{1}{2}+\dfrac{1}{2}i$,
we have
$$
p(z_0)
= \alpha+\frac{1-\alpha}{2}i.
$$

Putting $\beta = \dfrac{1-\alpha}{2}$,
we obtain
$$
\frac{z_0p'(z_0)}{p(z_0)}
= -\frac{4\alpha(1-\alpha)}{4\alpha^2+(1-\alpha)^2}+\frac{2(1-\alpha)^2}{4\alpha^2+(1-\alpha)^2}i
= -\frac{\alpha \beta k}{\alpha^2+\beta^2}+\frac{\beta^2 k}{\alpha^2+\beta^2}i
$$
where
$$
k
= 2
\geqq \frac{5}{4}
= \frac{1}{2}\left( \frac{\beta}{1-\alpha}+\frac{1-\alpha}{\beta} \right).
$$

For the case $z_0=-\dfrac{1}{2}-\dfrac{1}{2}i$,
we have
$$
p(z_0)
= \alpha-\frac{1-\alpha}{2}i.
$$

Putting $\beta = -\dfrac{1-\alpha}{2}$,
we obtain also
$$
\frac{z_0p'(z_0)}{p(z_0)}
= -\frac{4\alpha(1-\alpha)}{4\alpha^2+(1-\alpha)^2}-\frac{2(1-\alpha)^2}{4\alpha^2+(1-\alpha)^2}i
= -\frac{\alpha \beta k}{\alpha^2+\beta^2}+\frac{\beta^2 k}{\alpha^2+\beta^2}i
$$
where
$$
k
= -2
\leqq -\frac{5}{4}
= \frac{1}{2}\left( \frac{\beta}{1-\alpha}+\frac{1-\alpha}{\beta} \right).
$$

The function $p(z)$ satisfies Theorem \ref{p01thm1}.

Especially, the function
$$
p(z)
= 1+z+\frac{1}{2}z^2
\qquad(z\in\mathbb{U})
$$
is one of the example of Theorem Theorem \ref{p01thm1}.
In fact, when we choice a point $z_0$ such that
$$
z_0
= -\frac{1}{2}\pm\frac{1}{2}i
$$
for $|z_0|=\dfrac{1}{\sqrt[]{2}}$,
the function $p(z)$ satisfies that $\Re(p(z_0))>\dfrac{1}{2}$ for $|z|<|z_0|$ and $\Re(p(z_0))=\dfrac{1}{2}$.

For $z_0=-\dfrac{1}{2}+\dfrac{1}{2}i$,
we have
$$
p(z_0)
= \frac{1}{2}+\frac{1}{4}i
$$
and
$$
\frac{z_0p'(z_0)}{p(z_0)}
= -\frac{4}{5} +\frac{2}{5}i
= -\frac{2}{5}k +\frac{1}{5}ki
$$
with $k=2\geqq\dfrac{5}{4}$.
Furthermore, for $z_0=-\dfrac{1}{2}-\dfrac{1}{2}i$,
we get
$$
p(z_0)
= \frac{1}{2}-\frac{1}{4}i
$$
and
$$
\frac{z_0p'(z_0)}{p(z_0)}
= -\frac{4}{5} -\frac{2}{5}i
= \frac{2}{5}k +\frac{1}{5}ki
$$
with $k=-2\leqq-\dfrac{5}{4}$.

\

\begin{figure*}[h]
\label{p1fig1}
\begin{center}
\includegraphics[height=11.5cm]{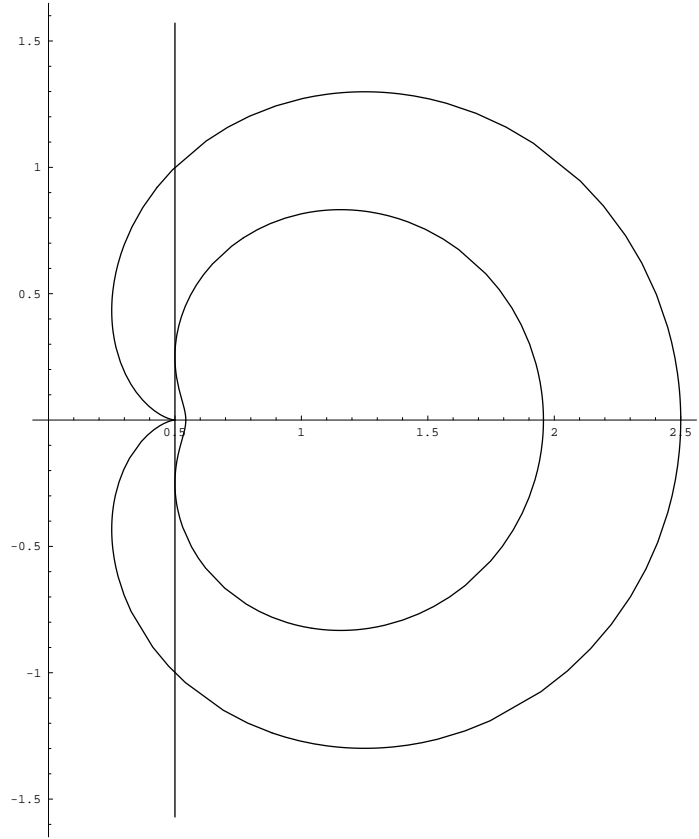}
\end{center}
\caption{$p(z)=1+z+\dfrac{1}{2}z^2$ in $|z|=1$ and $|z|=\dfrac{1}{\sqrt[]{2}}$}
\end{figure*}
\end{ex}

\

\end{document}